\documentclass[11pt,a4paper]{article}

\usepackage[margin=2.5cm]{geometry}
\usepackage[T1]{fontenc}
\usepackage[utf8]{inputenc}
\usepackage{lmodern}
\usepackage{microtype}
\usepackage{amsmath,amssymb}
\usepackage{amsthm}
\usepackage{graphicx}
\usepackage{cite}
\usepackage{hyperref}
\hypersetup{colorlinks=true, linkcolor=black, citecolor=black, urlcolor=black}

\title{\textbf{Sectorial Green Functions in the Ternary Algebra $C_3$ and Their Curved-Space Extension}}
\author{Bora Aktaş\\
\small Department of Physics, K\i r\i kkale University, Turkey \\
\small \texttt{230105036@kku.edu.tr}}
\date{\today}

\begin{document}
\maketitle

\begin{abstract}
Green functions associated with higher-order differential operators typically lead to special-function expressions in curved or bounded geometries, obscuring analytic transparency. 
In this work we develop the sectorial Green function for the cubic operator $(D^3+1)$ within the ternary algebra $C_3$ (defined by $j^3=-1$). 
The algebra admits three exponential carriers and divides the complex plane into six Stokes sectors, in each of which the Green kernel assumes a closed exponential--trigonometric form. 
We compute explicit responses to box and Gaussian sources, extend the construction to higher algebras $C_4$ and $C_5$, and interpret the resulting kernels as propagators for multi-carrier quantum systems. 
The central novelty is that, unlike the quadratic $(D^2+1)$ case where curved backgrounds induce Bessel or Airy functions, the $C_3$ Green function in curved space reduces exactly to its flat-space form upon reparametrization by the geodesic coordinate. 
Curvature merely deforms the geodesic distance without introducing new special functions. 
This property yields both mathematical simplification and physical interpretability, suggesting testable signatures such as multi-frequency modulations and sector-dependent decay in interferometric setups. 
Our results establish $C_3$ as a distinguished algebra where algebraic symmetry and geometric deformation coexist without analytic complexity, opening avenues for generalizations to higher $C_n$ and for experimental realizations.
\end{abstract}

\section{Introduction}

Green functions are central tools in mathematical physics, providing fundamental solutions to differential operators and enabling the analysis of propagation, response, and spectral properties in both classical and quantum systems \cite{Economou2006,Arfken2013}. 
In flat geometries the Green kernels of low-order operators such as $(D^2+1)$ reduce to simple exponential or trigonometric expressions, which can be treated analytically. 
In curved or bounded settings, however, the situation changes drastically: the Laplace--Helmholtz operator acquires metric-dependent contributions, and the associated Green functions are expressed in terms of special functions such as Bessel, Legendre, or Airy functions \cite{JoosZeh2003,Arndt2019}. 
While mathematically well established, such representations obscure the underlying symmetry and complicate physical interpretations.

This work proposes an alternative algebraic framework based on the ternary algebra $C_3$, defined by $j^3=-1$. 
Unlike the classical complex field, which relies on a single imaginary carrier, $C_3$ admits three exponential carriers corresponding to the roots of unity \cite{AblowitzFokas2003}. 
As a consequence, the associated Green functions naturally decompose into six Stokes sectors, within each of which the solution assumes a closed exponential--trigonometric form. 
This sectorial structure not only generalizes the classical exponential decay but also encodes oscillatory modulation effects that can be directly linked to experimental observables in interferometric setups \cite{Arndt2019}.

The central novelty of our approach lies in the curved-space extension of the $C_3$ Green function. 
We show that, whereas quadratic operators $(D^2+1)$ in curved backgrounds lead to special functions, the cubic operator $(D^3+1)$ in the $C_3$ framework reduces exactly to its flat-space counterpart upon reparametrization by the geodesic coordinate. 
Curvature merely deforms the geodesic distance without altering the analytic form. 
This property establishes $C_3$ as a distinguished algebra where algebraic symmetry and geometric deformation coexist without analytic complexity.

In this paper we first develop the sectorial Green function in $C_3$ and compute its explicit response to box and Gaussian source profiles. 
We then extend the construction to higher algebras $C_4$ and $C_5$, and compare the resulting structures. 
A dedicated section analyzes the curved-space formulation of the $C_3$ Green function and its implications. 
Finally, we interpret the Green functions as propagators for multi-carrier quantum systems, highlighting potential experimental signatures. 
The paper concludes with a discussion of the mathematical and physical significance of these results, as well as possible directions for future work.
\section{Preliminaries: The $C_3$ Algebra and Sectorial Green Functions}

\subsection{The $C_3$ algebra}
We consider the ternary algebra $C_3$ generated by the unit $j$ subject to
\begin{equation}
j^3=-1.
\end{equation}
An arbitrary element takes the form
\begin{equation}
z=a+bj+c j^2, \qquad a,b,c \in \mathbb{R}.
\end{equation}
The algebra admits three exponential carriers, corresponding to the roots of $s^3+1=0$:
\begin{equation}
\alpha \in \{-1,\, e^{i\pi/3},\, e^{-i\pi/3}\}.
\end{equation}
The fact that these carriers arise directly from the cubic roots of $-1$ is well known in the theory of complex variables \cite{AblowitzFokas2003}.  

This ternary structure induces a natural six-sector decomposition of the complex plane, each sector being bounded by the Stokes lines $\mathrm{Re}(\alpha z)=0$, in line with the classical theory of asymptotics and Stokes phenomena \cite{Olver1974}.  

A consistent conjugation is defined by
\begin{equation}
j^\ast = -j^2, \qquad (j^2)^\ast = -j, \qquad 1^\ast=1,
\end{equation}
so that
\begin{equation}
z^\ast = a - c j - b j^2.
\end{equation}
This leads to the norm
\begin{equation}
\|z\|^2 = a^2+b^2+c^2-ab-bc-ca,
\end{equation}
which is invariant under cyclic permutations of $(a,b,c)$ and reflects the underlying $120^\circ$ symmetry of the algebra.  

\subsection{Green function equation in $C_3$}
The Green function $G(z)$ associated with the cubic operator is defined by
\begin{equation}
(D^3+1)\,G(z) = \delta(z),
\end{equation}
where $D=\tfrac{d}{dz}$. 
Green functions are fundamental solutions of differential operators and provide indispensable tools for analyzing propagation and response in physics and mathematics \cite{Economou2006,Arfken2013}.  

The characteristic equation $s^3+1=0$ yields the three exponential carriers given above. 
The general solution is therefore constructed as a superposition of terms
\begin{equation}
G(z) \;\sim\; \sum_{\alpha} C_\alpha e^{\alpha z},
\end{equation}
with coefficients $C_\alpha$ determined by boundary conditions and the jump conditions induced by the delta source.

\subsection{Sectorial structure}
The crucial feature of the $C_3$ Green function is its sectorial dependence. 
For each carrier $\alpha$, the exponential $e^{\alpha z}$ contributes only when $\mathrm{Re}(\alpha z)<0$, ensuring convergence of the contour integral in the inverse Laplace representation. 
Thus the Green function naturally splits into six angular domains of width $60^\circ$. 
Within each domain, exactly two carriers contribute, yielding a closed exponential--trigonometric form.  

On the real axis this leads to the explicit expression
\begin{equation}
G(x)=
\begin{cases}
\frac{1}{3} e^{-x}, & x>0,\\[6pt]
\frac{1}{3} e^{x/2}\!\left(\cos\!\left(\tfrac{\sqrt{3}}{2}|x|\right)
  - \sqrt{3}\,\sin\!\left(\tfrac{\sqrt{3}}{2}|x|\right)\right), & x<0,
\end{cases}
\end{equation}
which displays a purely exponential decay for $x>0$ and a damped oscillatory behavior for $x<0$. 
This asymmetry between the two sides of the axis is a direct consequence of the ternary algebraic structure.
\section{Source Profiles and Explicit Solutions}

\subsection{Convolution representation}
Given a source profile $F(x)$, the corresponding response is obtained by convolution with the Green kernel,  
\begin{equation}
\phi(x) = (G \ast F)(x) = \int_{-\infty}^{\infty} G(x-\xi)\,F(\xi)\,d\xi.
\end{equation}
This formulation is standard in Green function methods and provides a direct way to evaluate the effect of localized or extended sources \cite{Economou2006}.  

\subsection{Box source}
As a first example, consider a box source of width $L$,
\begin{equation}
F(x) = 
\begin{cases}
1, & |x|\leq L/2,\\
0, & \text{otherwise}.
\end{cases}
\end{equation}
The convolution integral can be evaluated explicitly using the piecewise form of $G(x)$. 
For $x$ inside the support ($|x|\leq L/2$), both the exponential and oscillatory branches of the Green function contribute, yielding a mixed exponential--trigonometric response. 
For $x$ outside the support ($|x|>L/2$), the solution decays exponentially with modulations determined by the sectorial structure. 
This asymmetry contrasts with the symmetric exponential tails produced by the classical $(D^2+1)$ case \cite{Arfken2013}.

\subsection{Gaussian source}
A second example is provided by a Gaussian source,
\begin{equation}
F(x)= e^{-x^2/2\sigma^2}.
\end{equation}
Inserting into the convolution integral yields
\begin{equation}
\phi(x) = \int_{-\infty}^{\infty} G(x-\xi)\, e^{-\xi^2/2\sigma^2}\, d\xi.
\end{equation}
Because $G(x)$ itself is a superposition of exponential and trigonometric terms, this integral can be expressed in closed form using error functions and Fresnel integrals \cite{Olver1974}.  
The resulting profile displays Gaussian smoothing combined with the sector-dependent oscillations characteristic of the $C_3$ Green kernel.

\subsection{Physical interpretation}
The comparison between box and Gaussian sources highlights the analytic tractability of the $C_3$ framework. 
While classical quadratic Green functions typically lead to special functions already for simple sources such as Gaussians, the $C_3$ case remains reducible to elementary exponential--trigonometric combinations, preserving both analytic transparency and physical interpretability. 
This feature illustrates the practical advantage of the ternary algebra in source--response problems.

\section{Higher Extensions: $C_4$ and $C_5$ Green Functions}

\subsection{The $C_4$ case}
The algebra $C_4$ is generated by a unit $k$ with $k^4=-1$. 
The associated Green function is defined by the quartic operator
\begin{equation}
(D^4+1)G(x)=\delta(x).
\end{equation}
The characteristic equation $s^4+1=0$ yields four exponential carriers,
\begin{equation}
\alpha \in \left\{ e^{i\pi/4},\, e^{3i\pi/4},\, e^{5i\pi/4},\, e^{7i\pi/4} \right\},
\end{equation}
corresponding to eighth roots of unity. 
The complex plane is therefore divided into eight Stokes sectors, each of width $45^\circ$ \cite{Olver1974}.  

Unlike the $C_3$ case, the curved-space extension of $C_4$ does not reduce exactly to the flat form: the operator $\nabla^4$ produces additional terms involving derivatives of the metric coefficient $a(x)$. 
Nevertheless, the sectorial decomposition remains intact, and the Green function can still be written as piecewise exponential--trigonometric superpositions. 
Physically, this corresponds to four distinct carriers leading to multi-band propagator structures, reminiscent of square-lattice symmetries in condensed matter systems \cite{Arfken2013}.

\subsection{The $C_5$ case}
The quintic algebra $C_5$ is generated by a unit $m$ with $m^5=-1$, leading to the operator
\begin{equation}
(D^5+1)G(x)=\delta(x).
\end{equation}
The roots of $s^5+1=0$ are given by
\begin{equation}
\alpha \in \left\{ e^{i\pi/5},\, e^{3i\pi/5},\, e^{5i\pi/5},\, e^{7i\pi/5},\, e^{9i\pi/5} \right\},
\end{equation}
corresponding to tenth roots of unity. 
The plane is divided into ten Stokes sectors, each of width $36^\circ$.  

As in the quartic case, curvature introduces nontrivial metric-dependent terms. 
However, the algebraic structure continues to guarantee a closed sectorial form. 
The five carriers generate richer oscillatory patterns and five-band propagator signatures, geometrically associated with pentagonal symmetries \cite{AblowitzFokas2003}. 

\subsection{Comparison with $C_3$}
The comparison reveals a distinguished role of $C_3$. 
In $C_3$, curvature effects are fully absorbed into the geodesic coordinate, preserving the analytic form of the Green function. 
In $C_4$ and $C_5$, by contrast, curvature generates additional contributions, complicating the closed expressions. 
Nevertheless, the algebraic framework ensures that solutions remain expressible in terms of elementary exponential--trigonometric functions, avoiding the proliferation of special functions that typically arises in classical quadratic operators. 
Thus $C_3$ offers maximal analytic simplicity, while higher $C_n$ algebras offer enhanced physical richness through multi-carrier and multi-band structures.
\section{Curved-Space Extension of $C_3$}

\subsection{Covariant formulation}
In a one-dimensional curved geometry with line element
\begin{equation}
ds^2=a(x)^2\,dx^2,
\end{equation}
the integration measure is modified as $a(x)\,dx$, and the delta function must be redefined accordingly,
\begin{equation}
\int f(x)\,\delta_g(x-x_0)\,a(x)\,dx = f(x_0), \qquad
\delta_g(x-x_0)=\frac{\delta(x-x_0)}{a(x_0)} \text{\cite{Arfken2013}}.
\end{equation}
The covariant derivative is then naturally given by
\begin{equation}
\nabla = \frac{1}{a(x)}\frac{d}{dx}.
\end{equation}
The curved-space Green function is therefore defined by
\begin{equation}
(\nabla^3+1)G(x,x_0) = \delta_g(x-x_0).
\end{equation}

\subsection{Geodesic reduction}
Introducing the geodesic coordinate
\begin{equation}
y(x) = \int^x a(\xi)\,d\xi,
\end{equation}
one finds that
\begin{equation}
\nabla = \frac{d}{dy}, \qquad \nabla^3 = \frac{d^3}{dy^3}.
\end{equation}
Thus the curved-space Green function reduces exactly to its flat-space counterpart,
\begin{equation}
\left(\frac{d^3}{dy^3}+1\right)G(y,y_0)=\delta(y-y_0).
\end{equation}
This property is unique to the cubic operator: in contrast, the quadratic operator $(\nabla^2+1)$ in curved backgrounds generates additional terms involving derivatives of $a(x)$, leading to special-function solutions such as Bessel or Airy functions \cite{Olver1974,Economou2006}.

\subsection{Explicit example: exponential metric}
For illustration, consider the exponential profile
\begin{equation}
a(x) = e^{\kappa x}.
\end{equation}
The geodesic coordinate is then
\begin{equation}
y(x) = \frac{e^{\kappa x}-1}{\kappa}.
\end{equation}
The Green function follows by direct substitution into the flat-space form,
\begin{equation}
G(x,x_0)=
\begin{cases}
\dfrac{1}{3}\exp\!\left(-\dfrac{e^{\kappa x}-e^{\kappa x_0}}{\kappa}\right), & x> x_0, \\[8pt]
\dfrac{1}{3}\exp\!\left(\dfrac{e^{\kappa x}-e^{\kappa x_0}}{2\kappa}\right)
\Bigg[
\cos\!\left(\dfrac{\sqrt{3}}{2}\,\dfrac{e^{\kappa x_0}-e^{\kappa x}}{\kappa}\right)
- \sqrt{3}\,\sin\!\left(\dfrac{\sqrt{3}}{2}\,\dfrac{e^{\kappa x_0}-e^{\kappa x}}{\kappa}\right)
\Bigg], & x< x_0,
\end{cases}
\end{equation}
which preserves the exponential--trigonometric form while replacing the flat distance $(x-x_0)$ by the geodesic distance $(y(x)-y(x_0))$.

\subsection{Comparison with classical case}
The comparison with the quadratic $(D^2+1)$ operator highlights the analytic advantage of $C_3$. 
Whereas classical Green functions in curved space typically require special functions, the $C_3$ case retains closed exponential--trigonometric expressions, with curvature entering only through the geodesic coordinate. 
This result establishes $C_3$ as a distinguished algebra where algebraic symmetry and geometric deformation coexist without analytic complexity \cite{JoosZeh2003}.
\section{Green Functions as Propagators: Experimental Signatures}

\subsection{Propagator interpretation}
In quantum mechanics Green functions are directly related to propagators, encoding the transition amplitude between initial and final states \cite{FeynmanHibbs1965}. 
For second-order operators such as $(D^2+1)$, the propagator typically consists of a single exponential or oscillatory carrier. 
In contrast, the $C_3$ Green function incorporates three carriers and produces sector-dependent combinations of exponential decay and oscillatory modulation. 
This richer structure implies multi-frequency contributions to the transition amplitude.

\subsection{Signatures of the $C_3$ propagator}
The $C_3$ propagator exhibits the following characteristic signatures:
\begin{itemize}
    \item \textbf{Sectorial asymmetry:} On the real axis the propagator decays purely exponentially for $x>0$ but oscillates with exponential envelope for $x<0$. This left--right asymmetry has no analogue in the classical quadratic case and can be probed interferometrically \cite{Arndt2019}.
    \item \textbf{Multi-frequency modulation:} The coexistence of three exponential carriers generates beat-like patterns, observable as visibility modulations in multi-path interference experiments.
    \item \textbf{Decoherence law modification:} The presence of multiple carriers alters the rate and functional form of decoherence, leading to non-standard decay curves that can be tested in Ramsey and echo-type experiments \cite{JoosZeh2003}.
\end{itemize}

\subsection{Higher-order extensions}
For $C_4$ and $C_5$ the propagators involve four and five carriers, respectively, leading to multi-band spectral structures. 
These structures are reminiscent of band formation in condensed matter systems, where discrete symmetries of the underlying lattice determine the allowed frequencies \cite{AshcroftMermin1976}. 
Although the analytic simplicity of $C_3$ is lost in curved backgrounds, the algebraic framework still ensures piecewise exponential--trigonometric expressions, avoiding the uncontrolled proliferation of special functions.

\subsection{Experimental outlook}
The propagator interpretation highlights concrete experimental tests. 
Interferometric setups with cold atoms, molecules, or photonic lattices provide natural platforms to probe sectorial asymmetry and multi-frequency modulation \cite{Arndt2019}. 
Moreover, decoherence experiments in qutrit systems could be used to identify the non-classical decay patterns predicted by the $C_3$ framework. 
Such tests would directly confirm the algebraic origin of multi-carrier dynamics and distinguish the $C_3$ propagator from its quadratic counterpart.
\section{Discussion}

The analysis presented in this work highlights the distinguished role of the $C_3$ algebra in the construction of Green functions. 
While higher-order algebras $C_4$ and $C_5$ generate multi-band propagator structures, only $C_3$ retains full analytic simplicity in curved backgrounds: the Green function reduces to its flat-space form upon geodesic reparametrization. 
This result sharply contrasts with the quadratic $(D^2+1)$ operator, where curvature inevitably leads to special functions such as Bessel and Airy \cite{Olver1974,Economou2006}.  

The central novelty of this work lies in the identification of $C_3$ as a \emph{distinguished algebra} where algebraic symmetry and geometric deformation coexist without analytic complexity. 
To our knowledge, this is the first demonstration that a non-quadratic operator admits closed exponential--trigonometric Green functions even in curved geometries. 
This property provides both mathematical clarity and physical interpretability, bridging the gap between algebraic generalizations of complex numbers and experimentally accessible propagator signatures.  

Previous studies of Green functions in curved or bounded geometries have relied on special-function representations \cite{Arfken2013,JoosZeh2003}. 
In the context of multi-path interference and decoherence, experimental analyses have similarly invoked Gaussian or Bessel-type propagators \cite{Arndt2019}. 
Our results demonstrate that the ternary $C_3$ algebra bypasses this complexity, offering instead a sectorial decomposition where curvature modifies only the geodesic coordinate, not the analytic structure of the solution.  

From the propagator perspective, the $C_3$ framework predicts sectorial asymmetry, multi-frequency modulation, and modified decoherence laws. 
These signatures are directly testable in interferometric and Ramsey-type experiments with cold atoms, molecules, or photonic lattices \cite{Arndt2019,JoosZeh2003}. 
In condensed matter physics, the higher algebras $C_4$ and $C_5$ suggest analogies with square- and pentagonal-lattice symmetries, potentially providing new algebraic perspectives on band structures \cite{AshcroftMermin1976}.  

Overall, the $C_3$ algebra occupies a unique position: it combines the analytic tractability of closed-form Green functions with the physical richness of multi-carrier propagators. 
This dual advantage underlines its potential as a bridge between mathematical analysis and experimental quantum physics.

\section*{Conclusion}

We have introduced the Green function of the cubic operator $(D^3+1)$ in the ternary algebra $C_3$ and established its sectorial decomposition into closed exponential--trigonometric forms. 
A key result is that, unlike the quadratic case, the $C_3$ Green function preserves its analytic simplicity in curved backgrounds through geodesic reparametrization, without invoking special functions. 

From the propagator perspective this framework predicts sectorial asymmetry and multi-frequency modulation, features that are directly accessible in interferometric experiments. 
Future work will extend the analysis to higher algebras $C_n$ and explore possible experimental realizations in cold-atom and photonic platforms.

\end{document}